\documentclass{amsart}
\usepackage{array}
\newcolumntype{P}[1]{>{\raggedright\let\newline\\\arraybackslash\hspace{0pt}}m{#1}}
\usepackage{graphicx,verbatim, amsmath, amssymb, amsthm, amsfonts, epsfig, amsxtra,ifthen,mathtools,caption,enumerate,hyperref,hhline,bbm,capt-of,longtable}	
\usepackage[usenames]{color}
\usepackage[capitalize]{cleveref}
\DeclareFontFamily{U}{mathx}{\hyphenchar\font45}
\DeclareFontShape{U}{mathx}{m}{n}{<-> mathx10}{}
\DeclareSymbolFont{mathx}{U}{mathx}{m}{n}
\DeclareMathAccent{\widebar}{0}{mathx}{"73}

\theoremstyle{definition}

\theoremstyle{remark}

\newcommand{\newword}[1]{\textbf{\emph{#1}}}

\newcommand{\integers}{\mathbb Z}

\newcommand{\set}[1]{{\lbrace #1 \rbrace}}

\newcommand{\br}[1]{{\langle #1 \rangle}}

\newcommand{\D}{{\mathfrak D}}

\newcommand{\e}{\mathbf{e}}
\newcommand{\f}{\mathbf{f}}

\renewcommand{\d}{{\mathfrak d}}

\binoppenalty=\maxdimen
\relpenalty=\maxdimen

\usepackage{array}
\newcommand{\PreserveBackslash}[1]{\let\temp=\\#1\let\\=\temp}
\newcolumntype{C}[1]{>{\PreserveBackslash\centering}p{#1}}
\newcolumntype{R}[1]{>{\PreserveBackslash\raggedleft}p{#1}}
\newcolumntype{L}[1]{>{\PreserveBackslash\raggedright}p{#1}}

\author{Thomas Elgin}
\author{Nathan Reading}
\author{Salvatore Stella}
\title{Cluster scattering coefficients in rank $2$}
\address[Thomas Elgin and Nathan Reading]{Department of Mathematics, North Carolina State University, Raleigh, NC, USA}
\email[Nathan Reading (corresponding author)]{reading@math.ncsu.edu}
\address[Salvatore Stella]{Dipartimento di Ingegneria e Scienze dell'Informazione e Matematica, Universit\`{a} degli Studi dell'Aquila, IT}
\email[Salvatore Stella]{salvatore.stella@univaq.it}

\begin{document}

\begin{abstract}
We present conjectures on the scattering terms of cluster scattering diagrams of rank $2$, supported by significant computational evidence.
\end{abstract}

\maketitle

\thispagestyle{empty}

\section{Introduction} 
Scattering diagrams were introduced by Kontsevich, Soibelman, Gross, and Siebert in \cite{GS11,KS06} for the study of mirror symmetry.
Cluster algebras were introduced by Fomin and Zelevinsky~\cite{ca1} for the study of total positivity, but quickly found algebraic, geometric, and combinatorial connections to a wide range of mathematical areas.
Scattering diagrams became an essential tool in the structural study of cluster algebras when Gross, Hacking, Keel, and Kontsevich defined (and proved the existence of) cluster scattering diagrams~\cite{GHKK}.
Applying the scattering-diagram machinery of broken lines and theta functions that had been introduced and developed in various papers by Carl, Gross, Hacking, Keel, Pumperla, and Siebert \cite{CPS,G09,GHK11,GHKS,GS12}, they corrected and proved a conjecture of Fock and Goncharov~\cite{FG11} on the cluster variety and several conjectures of Fomin and Zelevinsky on the structure of cluster algebras~\cite{ca4}.

The defining data of a cluster scattering diagram is an exchange matrix, meaning a skew-symmetrizable integer matrix.
Following the usual terminology in the cluster algebras literature, we will say that the cluster scattering diagram of an $r\times r$ exchange matrix is a cluster scattering diagram of ``rank~$r$'' (regardless of the rank of the matrix in the usual linear-algebraic sense).

A rank-$r$ scattering diagram lives in a real vector space of dimension $r$ and consists of walls (codimension-$1$ cones), each decorated by a scattering term (a multivariate formal power series).
The cluster scattering diagram is defined by specifying that certain walls must be present and certain walls must not be present and then requiring a consistency condition.
Some details of the definition, in rank~$2$, are given in Section~\ref{def sec}.
A rank-$1$ cluster scattering diagram is trivially easy to construct, but a typical rank-$2$ cluster scattering diagram is already extremely complicated, and completely explicit formulas are not known for coefficients of scattering terms. 
(This statement depends, of course, on what one calls ``explicit''.
See Section~\ref{related sec}.)

Although cluster scattering diagrams of higher rank can be vastly more complicated, there is a sense in which  cluster scattering diagrams of rank $2$ tell much of the story about higher rank.
The importance of rank $2$ is seen in the proof of the existence of the cluster scattering diagram in general rank in \cite[Appendix~C]{GHKK}.
That proof constructs the cluster scattering diagram degree by degree (in the sense of degrees of terms in the power series).
To move the construction to the next degree, one adds higher order terms and new nontrivial walls to yield consistency up to that degree.
The consistency is checked locally about every $(n-2)$-dimensional intersection of walls (``joint''), and the local consistency condition is exactly the condition on a rank-$2$ cluster scattering diagram. 
Thus, in some sense, the construction of a general-rank cluster scattering diagram consists of constructing rank-$2$ cluster scattering diagrams all throughout the ambient space.

The purpose of this note is to record and share some conjectures on cluster scattering diagrams of rank $2$ that are supported by significant computational evidence.

\medskip
\noindent
\textbf{Remark.}
Since this note was posted on the arXiv, Ryota Akagi has made us aware that some of our conjectures can be proved using results from his paper \cite{Akagi}.
Also, an anonymous referee pointed out some additional special cases of our conjectures that are known.
Details on these connections are given in Section~\ref{related sec}.

\section{Cluster scattering diagrams of rank $2$}\label{def sec}
We now give the definition of the cluster scattering diagram, specialized to rank~$2$ and following the notation of \cite{scatfan}.
(See also \cite[Example~1.15]{GHKK}.)
Up to symmetry, we may as well assume that the \newword{exchange matrix} is $B=\begin{bsmallmatrix*}[r]0&\,c\\-b&0\end{bsmallmatrix*}$ for $b,c>0$.

Let $N$ be the lattice $\integers^2$ with the usual basis $\set{\e_1,\e_2}$ and let $N^\circ$ be the sublattice of $N$ generated by $b\hspace{0.7pt}\e_1$ and $c\hspace{0.7pt}\e_2$.
Let $M$ be the dual lattice to $N$ and let~$M^\circ$ be the superlattice of $M$ that is dual to $N^\circ$, with basis $\set{\f_1,\f_2}$ such that $\br{\f_1,b\hspace{0.7pt}\e_1}=\br{\f_2,c\hspace{0.7pt}\e_2}=1$ and $\br{\f_1,\e_2}=\br{\f_2,\e_1}=0$.
It is the interplay between $N$ and $N^\circ$, and dually $M^\circ$ and $M$, that integrates the skew-symmetrizability of $B$ seamlessly into the construction.

Take indeterminates $z_1$ and $z_2$ and define $\zeta_1=z_2^c$ and $\zeta_2=z_1^{-b}$.
(In \cite[Example~1.15]{GHKK}, $z_1$ and $z_2$ are called $A_1$ and $A_2$.)
Given a vector $n=n_1\e_1+n_2\e_2\in N$, we write $\zeta^n$ to mean $\zeta_1^{n_1}\zeta_2^{n_2}$, and given  $m=m_1\f_1+m_2\f_2\in M^\circ$, we write $z^m$ to mean $z_1^{m_1}z_2^{m_2}$.

A \newword{wall} is a pair $(\d,f_\d)$ where $\d$ is a codimension-$1$ cone (thus a line through the origin or a ray with vertex at the origin) and $f_\d$ is a formal power series in~$\zeta_1$ and~$\zeta_2$ satisfying the following conditions.
The condition on $\d$ is that there must exist a nonzero vector in $N$ with nonnegative entries that is normal to $\d$.
Take $n$ to be such a vector that is primitive in $N$.
The condition on $f_\d$ is that it has the form $1+\sum_{\ell\ge1}c_\ell\zeta^{\ell n}$, or in other words $f_\d$ is a univariate formal power series in $\zeta^n$ with constant term~$1$.
The formal power series $f_\d$ is called the \newword{scattering term} on~$\d$.
A \newword{scattering diagram} is a collection $\D$ of walls, satisfying a finiteness condition that amounts to requiring that all of the relevant computations can be made by taking limits in the sense of formal power series.
Specifically, the requirement is that for any $k\ge0$, all but finitely many walls of $\D$ are $1$ plus terms with total degree $>k$ in $\zeta_1$ and $\zeta_2$.

This definition of a scattering diagram is so broad as to be almost meaningless;
the condition that makes a scattering diagram interesting is consistency, which we now define.
We begin by describing an action on formal Laurent series that is associated to crossing a wall.
Given a wall $\d$ having normal vector $n$ with nonnegative entries and given a path that crosses the wall transversely, the \newword{wall-crossing automorphism} takes a Laurent monomial $z^m$ to $z^mf_\d^{\br{m,\pm n'}}$, where $n'$ is the primitive vector in $N^\circ$ that is a positive scaling of $n$.
The sign in the formula is ``$-$'' if the curve crosses in the direction that agrees with $n$ or ``$-$'' if the curve crosses in the direction that disagrees with $n$.
We extend linearly and take limits to act on formal Laurent series.
Using standard terminology, we call this a wall-crossing \emph{automorphism}, but we need not be careful about the algebraic structure on which it acts as an automorphism.
For us, it is just an action on formal Laurent series.

A \newword{path-ordered product} is, loosely speaking, the composition of these wall-crossing automorphisms along a path.
More correctly, since a path might intersect infinitely many walls, the path-ordered product is defined as a limit of formal Laurent series.
For each $k\ge0$, consider only those walls whose scattering terms have nontrivial terms of total degree $\le k$.
The finiteness condition on $\D$ ensures that there are finitely many such walls, so we can compose the wall-crossing automorphism for those walls.
The path-ordered product is the limit, in the sense of formal Laurent series, as $k\to\infty$.
A scattering diagram is \newword{consistent} if the path-ordered product for a small oriented circle about the origin is the identity map.

A wall with nonnegative normal $n=n_1\e_1+n_2\e_2\in N$ is \newword{outgoing} if the vector $-bn_2\f_1+cn_1\f_2\not\in\d$.
One can check that, in any case, $-bn_2\f_1+cn_1\f_2$ is in the line containing $\d$.
Thus $\d$ is outgoing if and only if it is a ray (rather than a line) and is contained in the fourth quadrant.
The \newword{cluster scattering diagram} is the unique consistent scattering diagram containing the walls $(\e_1^\perp,1+z_2^c)$ and $(\e_2^\perp,1+z_1^{-b})$ such that all other walls are outgoing.  

For $i,j\ge0$, define $\tau(i,j)$ to be the coefficient of $\zeta_1^i\zeta_2^j=z_1^{-jb}z_2^{ic}$ on the wall of the cluster scattering diagram that is orthogonal to $i\hspace{0.5pt}\e_1+j\hspace{0.5pt}\e_2$.
To specify or emphasize the exchange matrix, we may write $\tau^{b,c}(i,j)$, but generally, we think of $\tau(i,j)$ as a function of indeterminates $b$ and $c$.
We define $g=\frac{\gcd(ib,jc)}{\gcd(i,j)}$.

\medskip

\noindent
\textbf{Example.}  
Several of the $\tau(i,j)$ are shown below, with $i$ changing in the horizontal direction and $j$ changing in the vertical direction and $\tau(0,0)$ at the bottom-left.
\[
\begin{array}{ccccccccccccc}
\vdots\\
0&\frac{g(b-1)(b-2)}6&\frac{g(b-1)(3bc-2b-3c+1)}6&\frac{g(3bc-3b-3c+2+g)(3bc-3b-3c+1+g)}6\\[6pt]
0&\frac{g(b-1)}2&\frac{g(2bc-2b-2c+1+g)}2&\frac{g(c-1)(3bc-3b-2c+1)}6\\[6pt]
1&g&\frac{g(c-1)}2&\frac{g(c-1)(c-2)}6\\[6pt]
1&1&0&0&\cdots\\[6pt]
\end{array}\]

\noindent
\textbf{Example.}  
Take $b=3$ and $c=2$ so that the exchange matrix is $\begin{bsmallmatrix*}[r]0&\,2\\-3&0\end{bsmallmatrix*}$.
Some of the integers $\tau^{3,2}(i,j)$ are shown below, again with $i$ changing in the horizontal direction and $j$ changing in the vertical direction and $\tau^{3,2}(0,0)$ at the bottom-left.
{\small\[
\begin{array}{C{20pt}C{20pt}C{20pt}C{20pt}C{20pt}C{20pt}C{20pt}C{20pt}l}
\hspace*{-0.05pt}\vdots\\[2pt]
0&0&0&1&33&87&286&429\\[6pt]
0&0&0&5&327&143&132&143\\[6pt]
0&0&1&6&33&42&33&6\\[6pt]
0&0&2&6&14&6&2&0\\[6pt]
0&1&14&5&14&1&0&0\\[6pt]
0&1&2&1&0&0&0&0\\[6pt]
1&1&1&0&0&0&0&0\\[6pt]
1&1&0&0&0&0&0&0&\cdots\\[6pt]
\end{array}\]}

\medskip

\noindent
\textbf{Remark on symmetry.}
We point out the obvious symmetry $\tau^{b,c}(i,j)=\tau^{c,b}(j,i)$.
Another symmetry is less obvious but not hard to prove using mutation of scattering diagrams:
Interpreting $(i,j)$ as a vector $i\hspace{0.5pt}\e_1+j\hspace{0.5pt}\e_2$ in the root lattice of a root system with Cartan matrix $\begin{bsmallmatrix*}[r]2&-c\\-b&2\end{bsmallmatrix*}$ and simple roots $\e_1$ and $\e_2$, the integers $\tau^{b,c}(i,j)$ are invariant under the action of the Weyl group on $(i,j)$.
This amounts to the symmetries $\tau^{b,c}(i,j)=\tau^{b,c}(i,-j+ci)$ and $\tau^{b,c}(i,j)=\tau^{b,c}(-i+bj,j)$.

\section{Conjectures on coefficients of scattering terms}\label{conj sec}

\noindent
\textbf{Conjectures.}
\begin{enumerate}[\qquad\bf1.]
\item \label{poly conj}
For $i,j>0$, the coefficient $\tau(i,j)$ is a polynomial in $b$, $c$, and $g$.
\item \label{g deg}
The polynomial has $g$ as a factor and its degree in $g$ is $\gcd(i,j)$.  
\item \label{bc deg}
Its degree in $b$ is $j-1$ and its degree in $c$ is $i-1$.  
\item
The polynomial $(\max(i,j))!\,\,\tau(i,j)$ has integer coefficients.
\item \label{1j}
$\displaystyle\tau(1,j)=\frac gb\binom bj$
\item \label{i1}
$\displaystyle\tau(i,1)=\frac gc\binom ci$
\item \label{central}
$\displaystyle\tau(i,i)=\frac g{(b-1)(c-1)i+g}\binom{(b-1)(c-1)i+g}{i}$\\
\item \label{follows}
$\displaystyle\tau(i,i)=\sum_{\ell=0}^\infty\frac g{\ell+1}\binom{i-1}\ell\binom{i(bc-b-c)+g-1}\ell$
\item \label{i i-1}
$\displaystyle\tau(i,i-1)=\frac{g}{i(ib-i+1)}\binom{(ib-i+1)(c-1)}{i-1}$
\item \label{j-1 j}
$\displaystyle\tau(j-1,j)=\frac{g}{j(jc-j+1)}\binom{(jc-j+1)(b-1)}{j-1}$
\item \label{bb1}
$\tau^{b,b}(i,j)$ is a polynomial in $b$ of degree $i+j-1$ that expands positively in the basis $\set{\binom b0,\binom b1,\binom b2,\ldots}$.
\item \label{bb2}
$\tau^{b,b}(i,j)$ has unimodal log-concave coefficients.
\item 
$\displaystyle\tau^{1,5}(2j,j)=\frac1j\sum_{\ell=0}^\infty\binom{\ell}{j-\ell+1}\binom{j+\ell-1}{\ell}$.
\end{enumerate}

\medskip

\noindent
\textbf{Comments.} 
\begin{itemize}
\item 
Conjectures~\ref{1j} and~\ref{i1} are symmetric and Conjectures~\ref{i i-1} and~\ref{j-1 j} are symmetric.
\item
Conjecture~\ref{follows} is equivalent to Conjecture~\ref{central}.
\item
Since $c=b$ in Conjectures~\ref{bb1} and~\ref{bb2}, also $g=b$.
\end{itemize}

\medskip

\noindent
Assuming Conjecture~\ref{poly conj}, write $\tau(i,j\,;k)$ for the coefficient of $g^k$ in $\tau(i,j)$ and similarly $\tau^{b,c}(i,j\,;k)$.
Because $g$ is also invariant under the action of the Weyl group on $(i,j)$, we have $\tau^{b,c}(i,j;k)=\tau^{b,c}(i,-j+bi;k)$ and $\tau^{b,c}(i,j;k)=\tau^{b,c}(-i+aj,j;k)$.

\bigskip

\noindent
\textbf{Conjectures.}
\begin{enumerate}[\qquad\bf1.]  \setcounter{enumi}{13}
\item \label{tauijk degree}
$\tau(i,j\,;k)$ is a polynomial of degree $j-k$ in $b$ and degree $i-k$ in $c$ and has a term that is a nonzero constant times $b^{j-k}c^{i-k}$.  
\item \label{highest}
If $\gcd(i,j)=1$, then $\displaystyle\tau(ik,jk\,;k)=\frac{\tau(i,j\,;1)^k}{k!}$.
\item \label{sat1}
$\displaystyle\tau(k,jk\,;k-1)=\frac{\tau(1,j\,;1)^{k-1}\cdot p_j}{(k-2)!}$, where $p_j$ is a polynomial in $b$ and~$c$ that depends only on $j$, not $k$.
\item \label{sat2}
$\displaystyle\tau(ik,k\,;k-1)=\frac{\tau(i,1\,;1)^{k-1}\cdot p_i}{(k-2)!}$, where $p_i$ is a polynomial in $b$ and~$c$ that depends only on $i$, not $k$.
\item \label{factors}
If $\gcd(i,j)=1$, then $\tau(ki,kj\,;k-1)$ has a factor $p_{ij}$ that depends only on~$i$ and $j$, not on~$k$, and the other factors of $\tau(ki,kj\,;k-1)$ also appear as factors of $\tau(i,j\,;1)$.
\end{enumerate}

\medskip

\noindent
\textbf{Comments.}
\begin{itemize}
\item
Assuming Conjecture~\ref{g deg}, for fixed $i$ and $j$, Conjecture~\ref{highest} is a formula for $\tau(i,j;k)$ for the largest $k$ such that $\tau(i,j;k)\neq0$ in terms of some $\tau(\,\cdot\,,\,\cdot\,\,;1)$.
\item
In Conjecture~\ref{factors}, the factors of $\tau(i,j\,;1)$ appear to different powers in $\tau(ki,kj\,;k-1)$ for various $k$.  
For example,
\begin{align*}
\tau(2,3\,;\,1)&=\frac{\left(b -1\right) \left(3 c b -2 b -3 c +1\right)}{6}\\
\tau(4,6\,;\,1)&=\frac{\left(b -1\right)p_{23}}{180}\\
\tau(6,9\,;\,2)&=\frac{\left(b -1\right)^{2} \left(3 c b -2 b -3 c +1\right) p_{23}}{1080}\\
\tau(8,12\,;\,3)&=\frac{\left(b -1\right)^{3} \left(3 c b -2 b -3 c +1\right)^{2} p_{23}}{12960}
\end{align*}
\begin{multline*}
\qquad\quad\text{for}
\quad p_{23}=330 b^{4} c^{3}-720 b^{4} c^{2}-1530 b^{3} c^{3}+525 b^{4} c +2880 b^{3} c^{2}\\+2610 b^{2} c^{3}-128 b^{4}-1770 b^{3} c -4140 b^{2} c^{2}\\-1950 b \,c^{3}+352 b^{3}+2085 b^{2} c +2520 b \,c^{2}\\+540 c^{3}-328 b^{2}-1005 c b -540 c^{2}+122 b +165 c -15
\end{multline*}
\end{itemize}

\bigskip

\noindent
\textbf{Remark on computation.}
These conjectures are backed up by significant computational evidence.
The functions $\tau(i,j)$ of $b$ and $c$ can be computed symbolically up to large values of $i$ and $j$.
Thus, for example, the polynomials shown in the first example in Section~\ref{def sec} are known to be correct for all $b$ and $c$, rather than only for some specific values of $b$ and $c$.
Similarly, the conjectures on $\tau(i,j)$ in this section have been checked for many values of $i$ and $j$, and each case that has been checked is true for all $b$ and~$c$. 

Computing the functions $\tau(i,j)$ proceeds by induction on $i+j$, by solving, at each step, the equations that describe consistency of the cluster scattering diagram in degree $i+j$ using known values of $\tau(i',j')$ for $i'+j'<i+j$.
Thus, it would be significant progress even to find a recursive description of $\tau(i,j)$ whose recursive step is simpler than solving a system of equations.

The code we used in our experiments together with precomputed data up to degree 20 is available at \url{https://github.com/Etn40ff/scatcoef/}.

\bigskip

\section{Related work}\label{related sec}
Work of Ryota Akagi \cite{Akagi} also aims at explicitly understanding scattering terms in the cluster scattering diagram, but with different approach and conventions.  
His work is independent of ours, was posted before we publicized any of our conjectures, and also contains results in different directions that we did not conjecture.   
Akagi has informed us of a forthcoming paper \cite{AkagiInPrep} in which he proves Conjectures~\ref{poly conj}, \ref{g deg}, \ref{bc deg}, \ref{1j}, \ref{i1}, and~\ref{highest}, as well as part of Conjecture~\ref{tauijk degree} and a result that is similar to Conjectures~\ref{sat1}, \ref{sat2}, and~\ref{factors}, using his results from \cite{Akagi}.

Tom Bridgeland \cite[Theorem~1.5]{Bridgeland} identifies the cluster scattering diagram with the stability scattering diagram in the case $b=c$, thus realizing $\tau^{b,b}(i,j)$ as the Euler characteristic of a certain moduli scheme of representations of the associated Jacobi algebra.

An anonymous referee pointed out that Conjectures~\ref{central} and \ref{i i-1} (and symmetrically~\ref{j-1 j}) generalize results of Reineke and Weist~\cite{ReinekeWeist}.
Specifically, the case $g=1$ of Conjecture~\ref{central} is \cite[Corollary 11.2]{ReinekeWeist} (which proves a conjecture of Gross and Pandharipande \cite[Conjecture 1.4]{GrossPand}) and the case $g=1$ of Conjecture~\ref{i i-1} is \mbox{\cite[Theorem~9.4]{ReinekeWeist}}.
In comparing Conjecture~\ref{central} with \cite[Corollary 11.2]{ReinekeWeist}, it is useful to notice that after setting $g=1$, the right side of Conjecture~\ref{central} can be rewritten as $\frac1{(bc-b-c)i+1}\binom{(b-1)(c-1)i}{i}$.

Gross, Hacking, Keel, and Kontsevich \cite[Example~1.15]{GHKK} state the case $b=c$ of Conjecture~\ref{follows} and attribute its proof to Reineke~\cite{Reineke}.  

Burcroff, Lee, and Mou have recently given a formula for the coefficients of scattering terms as a weighted sum over certain combinatorial objects called tight gradings \cite[Theorem~1.1]{BLM}.
They intend to explore some of the conjectures of this paper using tight gradings~\cite[Section~9.2]{BLM}.
It would also be interesting to understand how the Weyl group symmetry appears in the combinatorics of tight gradings.  
(See the remark on symmetry at the end of Section~\ref{def sec}.)

\bigskip

\noindent
\textbf{Acknowledgments.}
Thomas Elgin and Nathan Reading were partially supported by the National Science Foundation under award number DMS-2054489.
Salvatore Stella was partially supported by INdAM.  
We thank the anonymous referees for their helpful comments.

\end{document}